\documentclass[a4paper,12pt]{amsart}
\usepackage{amssymb}
\usepackage{ifthen}
\usepackage{graphicx}

\newtheorem{thm}{Theorem}
\newtheorem{cor}{Corollary}
\newtheorem{lem}{Lemma}

\newtheorem{rem}{Remark}

\newtheorem{conj}{Conjecture}
\newtheorem{prob}{Problem}
\theoremstyle{definition}
\newtheorem{defn}{Definition}
\newtheorem{example}{Example}

\newcounter {own}
\def\theown {\thesection       .\arabic{own}}

\newenvironment{pf}[1][]{%
 \vskip 3mm
 \noindent
 \ifthenelse{\equal{#1}{}}%
  {{\slshape Proof. }}%
  {{\slshape #1.} }%
 }%
{\qed\bigskip}

\newcounter{alphabet}
\newcounter{tmp}
\newenvironment{Thm}[1][]{\refstepcounter{alphabet}%
\bigskip%
\noindent%
{\bf Theorem \Alph{alphabet}}%
\ifthenelse{\equal{#1}{}}{}{ (#1)}%
{\bf .} \itshape}{\vskip 8pt}

\makeatletter
\newcommand{\Ref}[1]{\@ifundefined{r@#1}{}{\setcounter{tmp}{\ref{#1}}\Alph{tmp}}}
\makeatother

\newenvironment{Lem}[1][]{\refstepcounter{alphabet}%
\bigskip%
\noindent%
{\bf Lemma \Alph{alphabet}}%
{\bf .} \itshape}{\vskip 8pt}

\newenvironment{Conje}[1][]{\refstepcounter{alphabet}%
\bigskip%
\noindent%
{\bf Conjecture \Alph{alphabet}}%
{\bf .} \itshape}{\vskip 8pt}

\newcommand{\es}{{\mathcal S}}

\newcommand{\IC}{{\mathbb C}}
\newcommand{\ID}{{\mathbb D}}

\newcommand{\CC}{{\mathcal C}}

\def\be{\begin{equation}}
\def\ee{\end{equation}}

\newcommand{\bee}{\begin{enumerate}}
\newcommand{\eee}{\end{enumerate}}

\newcommand{\blem}{\begin{lem}}
\newcommand{\elem}{\end{lem}}
\newcommand{\bthm}{\begin{thm}}
\newcommand{\ethm}{\end{thm}}
\newcommand{\bcor}{\begin{cor}}
\newcommand{\ecor}{\end{cor}}
\newcommand{\beg}{\begin{example}}
\newcommand{\eeg}{\end{example}}
\newcommand{\begs}{\begin{examples}}
\newcommand{\eegs}{\end{examples}}
\newcommand{\bdefe}{\begin{defn}}
\newcommand{\edefe}{\end{defn}}
\newcommand{\bprob}{\begin{prob}}
\newcommand{\eprob}{\end{prob}}
\newcommand{\bei}{\begin{itemize}}
\newcommand{\eei}{\end{itemize}}

\newcommand{\bcon}{\begin{conj}}
\newcommand{\econ}{\end{conj}}
\newcommand{\bcons}{\begin{conjs}}
\newcommand{\econs}{\end{conjs}}
\newcommand{\bprop}{\begin{propo}}
\newcommand{\eprop}{\end{propo}}
\newcommand{\br}{\begin{rem}}
\newcommand{\er}{\end{rem}}
\newcommand{\brs}{\begin{rems}}
\newcommand{\ers}{\end{rems}}
\newcommand{\bo}{\begin{obser}}
\newcommand{\eo}{\end{obser}}
\newcommand{\bos}{\begin{obsers}}
\newcommand{\eos}{\end{obsers}}
\newcommand{\bpf}{\begin{pf}}
\newcommand{\epf}{\end{pf}}
\newcommand{\ba}{\begin{array}}
\newcommand{\ea}{\end{array}}
\newcommand{\beq}{\begin{eqnarray}}
\newcommand{\beqq}{\begin{eqnarray*}}
\newcommand{\eeq}{\end{eqnarray}}
\newcommand{\eeqq}{\end{eqnarray*}}

\newcommand{\ra}{\rightarrow}

\newcommand{\ds}{\displaystyle}

\def\cc{\setcounter{equation}{0}   
\setcounter{figure}{0}\setcounter{table}{0}}

\begin{document}
\bibliographystyle{amsplain}
\title[Coefficient Conditions for Harmonic Univalent Mappings and Hypergeometric Mappings]{Coefficient Conditions for Harmonic
Univalent Mappings and Hypergeometric Mappings}

\author{S. V. Bharanedhar}
\address{S. V. Bharanedhar, Department of Mathematics,
Indian Institute of Technology Madras, Chennai--600 036, India.}
\email{
}

\author{S. Ponnusamy${}^{~\mathbf{*}}$}
\address{S. Ponnusamy, Department of Mathematics,
Indian Institute of Technology Madras, Chennai--600 036, India.}
\email{samy@iitm.ac.in}

\keywords{Harmonic, univalent, close-to-convex, starlike, and convex mappings,
coefficient estimates, Gaussian hypergeometric functions, Hadamard product (convolution).\\
${}^{\mathbf{*}}$ Corresponding author
}
\date{\today  
File: BPon1${}_{-}$final.tex}

\begin{abstract}
In this paper, we obtain coefficient criteria for a normalized harmonic function defined in the unit disk to be
close-to-convex and fully starlike, respectively. Using these coefficient conditions, we present different classes of
harmonic close-to-convex (resp. fully starlike) functions involving Gaussian
hypergeometric functions.
In addition, we present
a convolution characterization for a class of univalent harmonic functions discussed recently by Mocanu,
and later by Bshouty and Lyzzaik in 2010. Our approach provide examples of harmonic polynomials that
are close-to-convex and starlike, respectively.
\end{abstract}



\maketitle
\pagestyle{myheadings}
\markboth{S. V. Bharanedhar and S.Ponnusamy}{Harmonic Univalent Functions}
\cc

\section{Introduction and Two Lemmas}\label{BP1-sec1}
One of the basic coefficient inequalities states that if a normalized power series
$ f(z)=z+\sum_{n=2}^{\infty}a_nz^n$
satisfies the condition
\be\label{BP1Ieq1}
\sum_{n=2}^{\infty}n|a_n|\leq 1,
\ee
then $f$ is analytic in the unit disk $\mathbb{D}=\{z:\,|z|<1\}$
and ${\rm Re\,} f'(z)>0$ in $\mathbb{D}$ and hence the range
$f(\mathbb{D})$ is a close-to-convex domain.
We recall that a domain $D$ is
\textit{close-to-convex} if the complement of $D$ can be written as a union of non-intersecting
half-lines. Moreover, it is also well-known that each $f$ satisfying the
condition (\ref{BP1Ieq1}) implies that $|zf'(z)/f(z)-1|<1$ for $z\in \ID$ and, in particular,
$f\in \es^*,$ the class of starlike univalent functions
in $\ID$. One of the most natural questions is therefore to discuss its
analog coefficient conditions for complex-valued harmonic functions to be close-to-convex or starlike
in $\mathbb{D}.$

A complex-valued harmonic function $f=u+iv$ in $\mathbb{D}$
admits the decomposition $f=h+\overline{g},$ where both $h$ and $g$ are analytic in $\mathbb{D}$
(see \cite{Clunie-Small-84}). Here $h$ and $g$ are referred to as analytic
and co-analytic parts of $f$. A complex-valued harmonic function $z \mapsto f(z)=h(z)+\overline{g(z)}$
is \textit{locally univalent} if and only if the Jacobian $J_f$ is non-vanishing in $\ID$,
where $J_f(z)=|h'(z)|^2-|g'(z)|^2.$
For convenience, we let $f(0)=0$ and $f_z(0)=1$ so that every harmonic function $f$ in $\mathbb{D}$ can be
written as
\be\label{BP1Ieq2}
f(z)=z+\sum_{n=2}^{\infty}a_nz^n+\overline{\sum_{n=1}^{\infty}b_nz^n}:=h+\overline{g}.
\ee
We denote by ${\mathcal H}$ the class of all normalized harmonic functions $f$ in $\ID$ of this form.
The class of functions $f\in {\mathcal H}$ that are sense-preserving and univalent in $\mathbb{D}$ is denoted by $\es_H.$
Two interesting subsets of $\es_H$ are
$$ \es_H^0=\{f\in \es_H:\, b_1=f_{\overline{z}}(0)=0\} ~\mbox{ and }~\es =\{f\in \es_H:\, g(z)\equiv 0\}.
$$
In the recent years, properties of the class $\es_H$ together with its interesting geometric subclasses
have been the subject of investigations. We refer to the pioneering works of Clunie and Sheil-Small
\cite{Clunie-Small-84}, the book of Duren \cite{Duren:Harmonic} and the recent survey article of
Bshouty and Hengartner \cite{BsHen-04}. Let $\CC$, $\CC_H$, and $\CC_H^0$ denote the subclasses of
$\es$, $\es_H$, and $\es_H^0$, respectively, with close-to-convex images.
In \cite{Hiroshi-Samy-2010} the following result has been proved.

\begin{Lem} \label{BP1Ilem1}
Suppose that $f=h+\overline{g}$, where $h(z)=z+\sum _{n=2}^{\infty}a_nz^n$ and $g(z)=\sum _{n=1}^{\infty}b_nz^n$
in a neighborhood of the origin and $|b_1|<1$. If
\be\label{BP1Ieq3}
\sum _{n=2}^{\infty}n|a_n| + \sum _{n=1}^{\infty}n|b_n| \leq 1,
\ee
then $f\in {\mathcal C}_H^1$, where
${\mathcal C}_H^1=\{f \in {\mathcal S}_H:\, \mbox{${\rm Re\,} f_z(z)>|f_{\overline{z}}(z)|$ in $\mathbb D$}\}.
$
\end{Lem}

The condition (\ref{BP1Ieq3}) is easily seen to be sufficient  for $f\in {\mathcal C}_H^1$ if $a_n$ and $b_n$
are non-positive for all $n\geq 1$ ($a_1=1$). Since the proof is routine as in the analytic case, we omit its detail.

In \cite{Mocanu80} (see also \cite{Hiroshi-Samy-2010} for a slightly general result), Mocanu has shown that
functions in ${\mathcal C}_H^1$ are univalent in $\ID$. On the other hand, in \cite{Hiroshi-Samy-2010}
the authors have shown that each $f\in {\mathcal C}_H^1$
is indeed close-to-convex in $\mathbb D$.  In view of the information known for the class of
analytic functions, it is natural to ask whether the coefficient condition (\ref{BP1Ieq3})
is sufficient for $f$ to belong to $\es_H^*,$ where
$$\es_H^*=\{f\in \es_H:\, f(\ID) \mbox{ is a starlike domain with respect to the origin}\}.
$$
Functions in $\es_H^*$ are called starlike functions. In the sequel, we also need
$$\es_H^{*0}=\{f\in \es_H^*:\, f_{\overline{z}}(0)=0\}.
$$
Harmonic starlikness is not a hereditary property, because it is possible that for $f\in \es_H^*$,
$f(|z|<r)$ is not necessarily starlike for each $r<1$ (see \cite{Duren:Harmonic}).

\bdefe
A harmonic mapping $f\in {\mathcal H}$ is said
to be \textit{fully starlike} (resp. fully convex) if each $|z|<r$ is mapped onto a
starlike (resp. convex) domain (see \cite{CDO-CMFT-04}).
\edefe

Fully convex mappings are known to be fully starlike but not the converse as the function
$f(z)=z+(1/n)\overline{z}^n$ $(n\geq 2)$
shows. It is easy to see that the harmonic Koebe function $K$ with the dilation $\omega (z)=z$
is not fully starlike, although $K=H+\overline{G}\in {\mathcal S}_H^{*0}$, where
$$
H(z)=\frac{z-\frac{1}{2}z^2+\frac{1}{6}z^3}{(1-z)^3} ~\textrm{ and
}~ G(z)=\frac{\frac{1}{2}z^2+\frac{1}{6}z^3}{(1-z)^3}.
$$
For further details, we refer to \cite{CDO-CMFT-04}.

\bdefe 
We say that a continuously differentiable function $f$ in $\ID$
is {\it starlike} in $\ID$ if it is sense-preserving, $f(0)=0$, $f(z)\neq0$
for all $z\in \ID\setminus \{0\}$ and
$$
{\rm Re} \left ( \frac{Df(z)}{f(z)}\right ) >0
\quad \mbox{for all $z\in \ID\setminus \{0\}$},
$$
where $Df=zf_{z}-\overline zf_{\overline z}.$
\edefe

The last condition gives that  ($z=re^{i\theta}$)
$$\frac{\partial}{\partial \theta}\big(\arg f(re^{i\theta})\big)
= {\rm Re} \left ( \frac{Df(z)}{f(z)}\right )  >0
\quad \mbox{for all $z\in \ID\setminus \{0\}$},
$$
showing that the curve $C_r=\{f(re^{i\theta}):\, 0\leq \theta <2\pi\}$
is starlike with respect to the origin for each $r\in(0, 1)$
(see \cite[Theorem 1]{Mocanu80}). In this case,  the last condition implies that
$f$ is indeed fully starlike in $\ID$.
At this place, it is also important to observe that
$Dg$ for $C^1$-functions behaves much like $zg'$ for
analytic functions, for example in the sense that for $g$ univalent and analytic in $\ID$,
$g$ is starlike if and only if ${\rm Re}\, (zg'(z)/g(z))> 0$ in $\ID$.
A similar characterization has also been obtained by  Mocanu \cite{Mocanu80} for convex ($C^2$)
functions.
It is worth to point out that in the case of analytic functions fully starlike (resp. fully convex)
is same as starlike (resp. convex) in $\ID$.
Lately, interesting distortion theorems and coefficients
estimates for convex and close-to-convex harmonic mappings were given by
Clunie and Sheil-Small \cite{Clunie-Small-84}.

As a consequence of convolution theorem \cite[Theorem~ 2.6, p.~908]{AhJaSi-2003} (see also
the proof of Theorem 1 in \cite{Jay-99}) these authors obtained
a sufficient coefficient condition for harmonic starlike mappings. Unfortunately, there is
a minor error in the main theorem and would like to point this out as we use this
for our applications.

\blem\label{BP1-lem1}
Let $f=h+\overline{g}\in \es_H^0$. Then $f$ is fully starlike in $\ID$ if and only if
\be\label{BP1Ieq4}
h(z)*A(z) -\overline{g(z)}* \overline{B(z)} \neq 0 \quad \mbox{ for $|\zeta|=1$, $0<|z|<1$,}
\ee
where
$$A(z)=\frac{z+((\zeta-1)/2)z^2}{(1-z)^2} ~\mbox{ and }~ B(z)=\frac{\overline{\zeta}
z-((\overline{\zeta}-1)/2)z^2}{(1-z)^2}.
$$
\elem
\bpf
A necessary and sufficient condition for a function $f$ to be starlike in $|z|<r$ for each $r<1$ is that
\be\label{BP1Ieq5}
\frac{\partial}{\partial \theta}(\arg f(re^{i\theta}))  =
{\rm Re\,} \left(\frac{zh'(z)-\overline{zg'(z)}}{h(z)+\overline{g(z)}}\right)>0
\quad \mbox{for all $z\in \ID\setminus \{0\}$}.
\ee
We remind the reader that if $f=h+\overline{g}\in \es_H$ with $g'(0)=b_1\ne 0$, then the limit
$$\lim_{z\ra 0} \frac{zh'(z)-\overline{zg'(z)}}{h(z)+\overline{g(z)}}
$$
does not exist, but the limit does exist which is $1$ when $b_1= 0$. This observation is crucial in
the remaining part of our proof. Thus, by (\ref{BP1Ieq5}), $f$ is fully starlike in $\ID$ if and only if
$$\frac{zh'(z)-\overline{zg'(z)}}{h(z)+\overline{g(z)}} \neq \frac{\zeta -1}{\zeta +1},
~~|\zeta |=1, ~\zeta \neq -1, ~0<|z|<1
$$
and, as in the proof of  Theorem 2.6 \cite{AhJaSi-2003}, a simple computation shows that the last
condition is equivalent to (\ref{BP1Ieq4}). The proof is complete.
\epf

In view of Lemma \ref{BP1Ilem1}, the hypothesis that $f=h+\overline{g}\in \es_H$ in
\cite[Corollary 2.7, p.~908]{AhJaSi-2003} can be relaxed as the condition (\ref{BP1Ieq3})
implies that $f\in \es_H.$
So we may now reformulate  \cite[Corollary 2.7, p.~908]{AhJaSi-2003}
in the following improved form (see also \cite{Silverman-98}).

\blem\label{BPlIlem2}
Let $f=h+\overline{g}$ be a harmonic function of the form {\rm (\ref{BP1Ieq2})}
with $b_1=g'(0)=0$. If
\be\label{BP1Ieq6}
\ds \sum _{n=2}^{\infty}n|a_n| + \sum _{n=2}^{\infty}n|b_n| \leq 1
\ee
then $f\in  {\mathcal C}_H^1\cap \es_H^{*0}$. Moreover, $f$ is fully starlike in $\ID$.
\elem
\bpf
By Lemma \Ref{BP1Ilem1}, the coefficient condition (\ref{BP1Ieq6}) ensures the univalency of $f$
and moreover, $f\in  {\mathcal C}_H^1$. Now, in order to show that
(\ref{BP1Ieq6}) implies $f\in \es_H^{*0}$, we apply Lemma \ref{BP1-lem1}. As in the proof
of \cite[Corollary~2.7]{AhJaSi-2003}, it suffices to show that the condition (\ref{BP1Ieq4}) holds.
Indeed, we easily have

\vspace{6pt}

$\left |h(z)*A(z) -\overline{g(z)}* \overline{B(z)}\right |$
\beqq
&=& \left | z+\sum_{n=2}^{\infty}\left (n+\frac{(n-1)(\zeta -1)}{2}\right )a_nz^n-
\sum_{n=2}^{\infty}\left (n\zeta-\frac{(n-1)(\zeta -1)}{2}\right )\overline{b_n}\overline{z^n}
\right |\\
&> & |z|\,\left [ 1 -\sum_{n=2}^{\infty}n|a_n| -\sum_{n=2}^{\infty}n|b_n|
\right ] \geq 0
\eeqq
and so Lemma \ref{BP1-lem1} gives that $f$ is fully starlike in $\ID$ and hence, $f\in \es_H^{*0}$.
\epf

For instance, according to Lemma \ref{BPlIlem2}, it follows that
if $\alpha\in \IC$ is such that $|\alpha|\leq 1/n$ for some $n\geq 2$, then the
function $f$ defined by
$$f(z)=z +\alpha \overline{z}^n
$$
belongs to ${\mathcal C}_H^1\cap \es_H^{*0}$.
Later in Section \ref{BP1-sec3}, we present a number of interesting
applications of Lemma \ref{BPlIlem2}.

\section{Conjecture of Mocanu on Harmonic Mappings}\label{BP1-sec2}

According to our notation, the conjecture of Mocanu  \cite{Mocanu-2010}
may be reformulated in the following form.

\begin{Conje}
If
$${\mathcal M} =\left \{f=h+\overline{g}\in {\mathcal H}:\,  ~g' =zh',~
{\rm Re\,}\left(1+z\frac{h''(z)}{h'(z)}\right)>-\frac{1}{2} ~\mbox{for $z\in \ID$}
\right\},
$$
then  $f\in \es_H^0$.
\end{Conje}

In \cite{Bshouty-Lyzzaik-2010}, Bshouty and Lyzzaik have solved the conjecture of Mocanu by
establishing the following stronger result.

\begin{Thm}\label{BP1Ath2}
${\mathcal M}\subset {\mathcal C}_H^0$.
\end{Thm}

It is worth to reformulate this result in a general form.

\begin{thm}\label{BP1Cth3}
Let $f=h+\overline{g}$ be a harmonic mapping of $\ID,$ with $h'(0)\neq0,$ that satisfies
$$ g'(z)=e^{i\theta}zh'(z) ~\mbox{ and }~
{\rm Re}\left(1+z\frac{h''(z)}{h'(z)}\right)>-\frac{1}{2} ~\mbox{ for all $z\in \ID$}.
$$
Then $f$ is a univalent close-to-convex mapping in $\ID$.
\end{thm}
\bpf
This theorem is proved for $\theta=0$ by Bshouty and Lyzzaik \cite{Bshouty-Lyzzaik-2010},
i.e. ${\mathcal M}\subset {\mathcal C}_H^0$.
However, it can be easily seen from their proof that the theorem continues to hold if the
dilatation $\omega $ is chosen to be $\omega(z)=e^{i\theta}z$ instead of $\omega(z)=z$.
So we omit the details.
\epf

Since the function $f\in {\mathcal M}$ satisfies the condition $f_{\overline{z}}(0)=0$, it is natural
to ask whether ${\mathcal M}$ is included in $\es_H^{*0}$ or in ${\mathcal C}_H^1$. First we construct
a function $f\in{\mathcal M}$ such that $f\not\in {\mathcal C}_H^1$.

\begin{figure}
\begin{center}
\includegraphics[width=9cm]{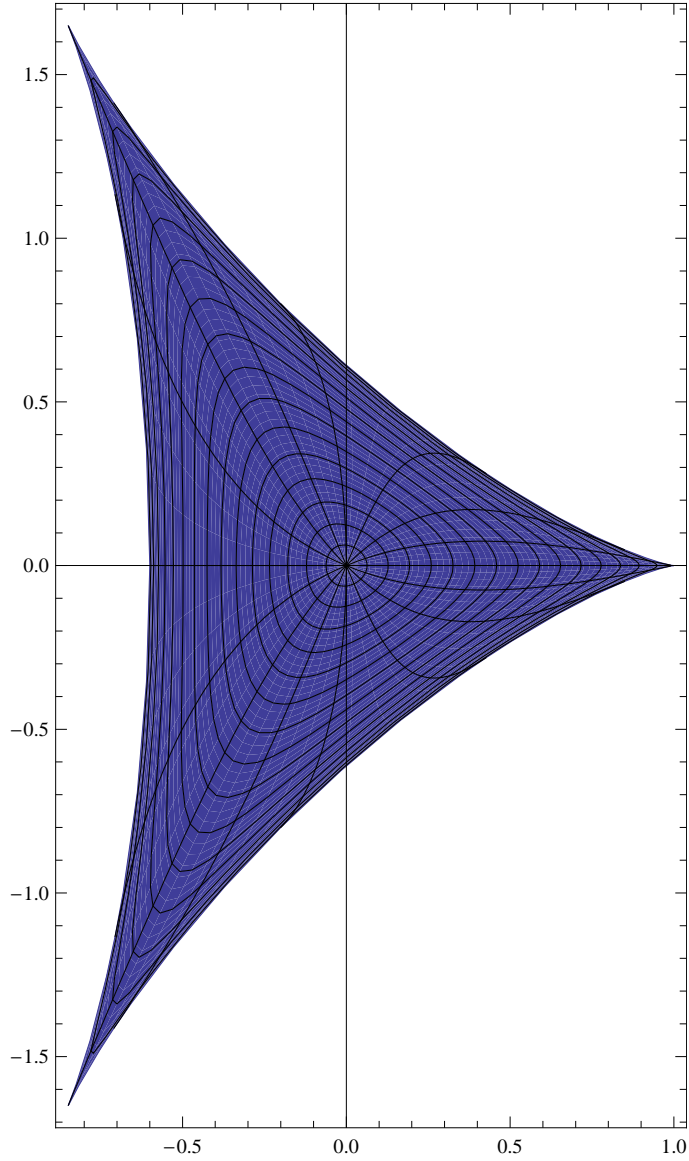}
\caption{The graph of the function $f_0(z)
=z-\frac{3}{10}z^2+\overline{\frac{z^2}{2}-\frac{1}{5}z^3}$\label{BP1-fig1}}
\end{center}
\end{figure}

Consider $f=h+\overline{g}$, where
$$h(z)=z-az^n  ~\mbox{ and }~ g(z)=\frac{z^2}{2}-\frac{n}{n+1}az^{n+1}
$$
for $n\geq 2$ and $0<a\leq1/n$.  It follows that $g'(z)=zh'(z)$ and
$$1+z\frac{h''(z)}{h'(z)}=\frac{1-n^2az^{n-1}}{1-naz^{n-1}}.
$$
Also it is a simple exercise to see that
$$w=\frac{1-n^2az^{n-1}}{1-naz^{n-1}}
$$
maps the unit disk $\ID$ onto the disk
$$\left|w-\frac{1-n^3a^2}{1-n^2a^2}\right|<\frac{an(n-1)}{1-n^2a^2}
$$
if $0<a <1/n$, and onto the half-plane ${\rm Re}\,w<(n+1)/2$ if $a=1/n$.
In particular this disk lies in the half-plane
$${\rm Re\,} w>\frac{1-n^2a}{1-na}
$$
and thus, ${\rm Re\,} w>-1/2$ if $(1-n^2a)/(1-na)\geq-1/2$, i.e. if $0<a\leq3/(n(1+2n))$.
According to Theorem \Ref{BP1Ath2}, $f=h+\overline{g}$ is univalent close-to-convex mapping
in $\ID$ whenever $a$ satisfies the condition
$$ 0<a\leq \frac{3}{n(1+2n)}.
$$
On the other hand, this function does not satisfy the coefficient condition
(\ref{BP1Ieq6}). Moreover, it can be easily seen that
$f\notin {\mathcal C}_H^1$. Indeed, if $a=0.3$ and $n=2,$ then the corresponding function
$$f_0(z)=z-\frac{3}{10}z^2+\overline{\frac{z^2}{2}-\frac{1}{5}z^3}
$$
does not belong to ${\mathcal C}_H^1$. The graph of $f_0(z)$ is shown in Figure \ref{BP1-fig1}.
This example shows that there are functions in ${\mathcal M}$ but are not necessarily
belonging to ${\mathcal C}_H^1$. Indeed the above discussion gives

\bthm
${\mathcal M}\not\subset {\mathcal C}_H^1$.
\ethm

Moreover the graph of
$$\ds f(z)=z-\frac{3}{n(2n+1)}z^n+\overline{\frac{z^2}{2}-\frac{3}{(n+1)(2n+1)}z^{n+1}},
$$
for various values of $n\geq2$, shows that $f(z)$ is starlike in $\ID$. This motivates us to state
\bcon
${\mathcal M}\subset {\es}_H^{*0}$.
\econ
Our next result gives a convolution characterization for functions $f\in {\mathcal M}$ to be
starlike in $\ID$.

\bthm 
Let $f=h+\overline{g}\in \es_H^0$ such that $g'(z)=zh'(z)$. Then 
$f$ is fully starlike in $\ID$ if and only if
\be\label{BP1Aeq1}
h(z)*A(z) -\overline{z}\left( \overline{h(z)}*\overline{B(z)}\right)\neq 0
\quad \mbox{ for $|\zeta|=1$, $0<|z|<1$},
\ee
where
$$ A(z)= \frac{2z+(\zeta-1)z^2}{(1-z)^2}  ~\mbox{ and }~
B(z)= \frac{2z^2+z(\overline{\zeta}-1)+(1-z)^2(\overline{\zeta}-1)\log(1-z)}{z(1-z)^2}.
$$
\ethm
\bpf As in the proof of Lemma \ref{BP1-lem1}, $f$ is fully starlike if and only if
\be\label{BP1Aeq2}
{\rm Re\,} \left(\frac{zh'(z)-\overline{zg'(z)}}{h(z)+\overline{g(z)}}\right)>0
\quad \mbox{for all $z\in \ID\setminus \{0\}$}.
\ee
Since $g'(0)=0$ and $g'(z)=zh'(z)$, we obtain  that
$$\lim_{z\ra 0} \frac{zh'(z)-\overline{zg'(z)}}{h(z)+\overline{g(z)}}=1,
$$
and therefore the condition (\ref{BP1Aeq2}) holds if and only if
$$\frac{zh'(z)-\overline{z^2h'(z)}}{h(z)+\overline{\int_0^zth'(t)\,dt}}
\neq \frac{\zeta -1}{\zeta +1} ~\mbox{ for $|\zeta |=1$, $\zeta \neq -1$, $0<|z|<1$}.
$$
The last condition is equivalent to
$$0 \neq(\zeta+1)\left[zh'(z)-\overline{z^2h'(z)}\right]-(\zeta-1)\left[h(z)+\overline{\int_0^zth'(t)\,dt}\right]
$$
which is same as
\be\label{BP1Aeq3}
0\neq h(z)* A(z)-\overline{\left[(\overline{\zeta}+1)z^2h'(z)+(\overline{\zeta}-1)\int_0^zth'(t)\,dt\right]}.
\ee
Finally, as
$$ z^2h'(z)=z\left[h(z)*\frac{z}{(1-z)^2}\right]$$
and
$$g(z)=\int_0^zth'(t)\,dt =\frac{z}{2}\left[h(z)*\left(\frac{2}{1-z}+\frac{2}{z}\log(1-z)\right)\right],
$$
the condition (\ref{BP1Aeq3}) is easily seen to be equivalent to the required convolution
condition (\ref{BP1Aeq1}). The proof is complete.
\epf

\begin{figure}
\begin{center}
\includegraphics[width=9cm]{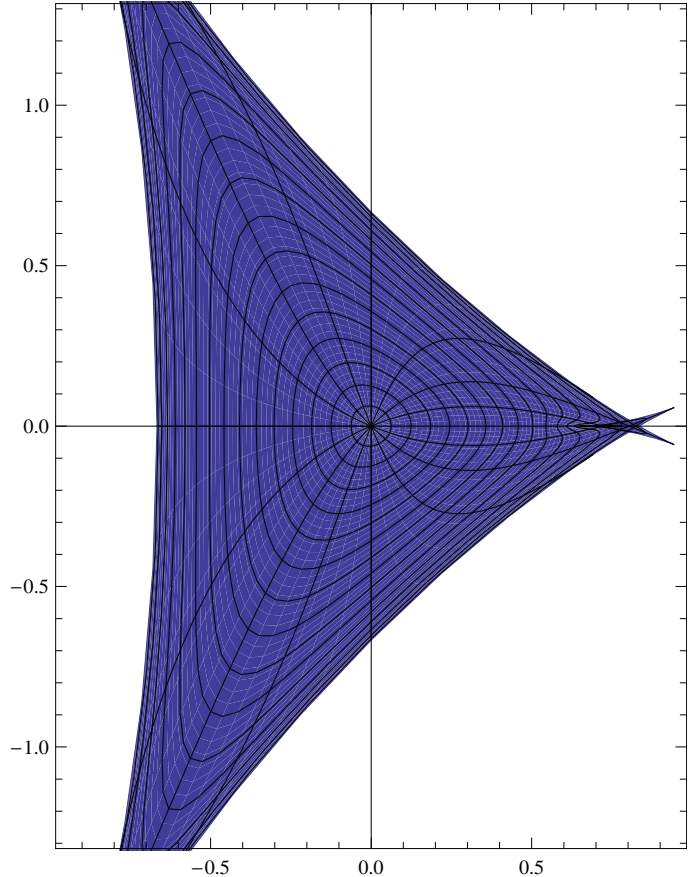}
\caption{Graph of the function $\ds f(z)=z- (1/2)z^2 +\overline{(1/2)z^2-(1/3)z^3}$\label{BP1-fig2}}
\end{center}
\end{figure}

Now, we consider the harmonic function $f=h+\overline{g}$, where
$$h(z)=z-\frac{z^2}{2} ~\mbox{ and  }~ g(z)=\frac{z^2}{2}-\frac{z^3}{3}
$$
so that $g'(z)=zh'(z)$. It follows that
$$1+z\frac{h''(z)}{h'(z)}=\frac{1-2z}{1-z}
$$
and we see that
\be\label{BP1Aeq5}
{\rm Re\,}\left(1+z\frac{h''(z)}{h'(z)}\right)<\frac{3}{2} ~\mbox{ for $z\in \ID$}.
\ee
The function $h$ satisfying the condition (\ref{BP1Aeq5}) is known to
satisfy the condition (see eg. \cite{PonRaj95,PoSinh96})
$$\left | \frac{zh'(z)}{h(z)}-\frac{2}{3} \right |<\frac{2}{3}, \ z\in \ID
$$
and hence, $h$ is starlike in $\ID$. The graph of $f(z)=h(z)+\overline{g(z)}$ shown in
Figure \ref{BP1-fig2} shows that $f=h+\overline{g}$ is not univalent in $\ID$.
This example motivates to raise the following

\bprob
For $\alpha \in (2/3,1)$, define
$${\mathcal P}(\alpha) =\left \{f=h+\overline{g}\in {\mathcal H}:\,  ~g' =zh',~
{\rm Re\,}\left(1+z\frac{h''(z)}{h'(z)}\right)<\frac{3\alpha}{2} ~\mbox{for $z\in \ID$}
\right\}.
$$
Determine $\inf\{\alpha \in (2/3,1):\, {\mathcal P}(\alpha) \subset \es_H^0\}$.
\eprob

\section{Harmonic Polynomials}\label{sec-harmpol}
One of the interesting problems in the class of harmonic mappings is to find a method of
constructing sense-preserving harmonic polynomials that have some interesting geometric properties.
In \cite{Suffridge-1998, Suffridge-2001}, the authors
discussed such polynomials with many interesting special cases. Prior to the work of
Suffridge \cite{Suffridge-1998}, only fewer examples of such polynomials were known.
In this section, we shall see that some
of the results of \cite{Suffridge-1998, Suffridge-2001} have closer link with our results in
Section \ref{BP1-sec1}. Following the ideas from \cite{Suffridge-1998, Suffridge-2001}, let
$Q(z)=\sum_{k=1}^nc_kz^k$ be a polynomial of degree $n$.
Define
$$\hat{Q}(z)=z^n\overline{Q(1/\overline{z})}. 
$$
Thus, if $Q(z)=c\prod_{j=1}^{n}(z-z_j)$ then
$\hat{Q}(z)=\overline{c}\prod_{j=1}^{n}(1-z\overline{z_j})$ and it follows that the zeros of $Q$
and $\hat{Q}$ on the unit circle $|z|=1$ are same. In \cite{Suffridge-1998} Suffridge proved the
following theorem.

\begin{Thm}{\rm (\cite[Theorem 1]{Suffridge-1998})}\label{BP1Cth1}
Let $Q(z)$ be a polynomial of degree $q\leq n-2$ with $Q(0)=1$ and assume that $Q(z)\neq0$ when $z\in\ID.$  Let
$h$ and $g$ be defined by $h(0)=g(0)=0$,  and
$$h'(z)=Q(z)+e^{i\phi}(1-t)z\hat{Q}(z),~  g'(z)=e^{i\beta}tz\hat{Q}(z)
$$
where $\phi$, $\beta$ and $t$ are real, $0\leq t\leq 1.$ Then the harmonic polynomial $f=h+\overline{g}$
has degree $n$ and is sense-preserving in $\ID$.
\end{Thm}

With an additional condition on $Q,$ we can improve this result by showing that the harmonic polynomial
$f=h+\overline{g}$ described in Theorem \Ref{BP1Cth1} is indeed close-to-convex in $\ID$.
More precisely we prove the following theorem.

\begin{thm}\label{BP1Cth2}
Let $Q$, $h$, $g$, $\phi$, $\beta$,  and $t$ be defined as in Theorem {\rm \Ref{BP1Cth1}}.
If  $Q$
satisfies the condition ${\rm Re\,}\{Q(z)\}>|z \hat{Q}(z)|$ for all
$z\in\ID$ then $f=h+\overline{g}$ belongs to $\mathcal{C}_H^1$ and hence, $f$ is close-to-convex in $\ID$.
\end{thm}
\bpf
It follows from the hypotheses that
\begin{align*}
 {\rm Re\,}(h'(z))&={\rm Re\,}(Q(z)+e^{i\phi}(1-t)z\hat{Q}(z))\\
&={\rm Re\,}(Q(z))+{\rm Re\,}(e^{i\phi}(1-t)z\hat{Q}(z))\\
&\geq {\rm Re\,}(Q(z))-|e^{i\phi}(1-t)z\hat{Q}(z)|\\
&>|z \hat{Q}(z)|-(1-t)|z \hat{Q}(z)|=|g'(z)|.
\end{align*}
Thus desired conclusion follows (see \cite{Hiroshi-Samy-2010}).
\epf

\beg\label{BP1Ceg1}
Consider
$$f(z)=z+e^{i\phi}\frac{(1-t)}{n}z^n+e^{i\beta}\frac{t}{m}\overline{z}^m, ~\mbox{ $n\geq2$, $m\geq1$,
$\phi\in\mathbb{R}$ and $\beta\in\mathbb{R}$},
$$
$0\leq t\leq 1.$ Then according to  Lemma \ref{BP1Ilem1}, we have
$$n|a_n|+m|b_m|=n\left|e^{i\phi}\frac{(1-t)}{n}\right|+m\left|e^{i\beta}\frac{t}{m}\right|=(1-t)+t=1
$$
showing that $f$ is not only close-to-convex, but also in ${\mathcal C}_H^1$.
On the other hand by Lemma \ref{BPlIlem2}, $f$ is also fully starlike in $\ID$ whenever $m\geq2$.
In particular the function
$$f(z)=z+e^{i\phi}(1-t)(z^2/2)+e^{i\beta}t(\overline{z}^2/2)
$$
is close-to-convex and fully starlike in $\ID$. By a direct method, Suffridge \cite[Example 1]{Suffridge-1998}
showed that this function is univalent in $\ID$.
\eeg

%

Using Theorem \ref{BP1Cth3}, it is possible to give a new proof of limit mapping theorem of
Suffridge et. al.  \cite[Theorem 3.1]{Suffridge-2001}. To do this, we assume that all the zeros of $Q(z)$
lie on the unit circle $|z|=1$. Then, for $q=n-2$ and $t=1$ in Theorem \Ref{BP1Cth1},
we have
$$ h'_n(z)=Q(z)=\prod_{j=1}^{n-2}(1-e^{-i\psi_j}z)=\frac{1-z^{n+1}}{\prod_{j=1}^{3}(1-ze^{i\psi_j})}
~\mbox{ and }~ g'_n(z)=z\prod_{j=1}^{n-2}(z-e^{i\psi_j}).
$$
It is clear that $h'_n(z)$ converges uniformly on the compact subsets of the unit disk to
\beq\label{BP1Ceq2}
h'(z)=\frac{1}{\prod_{j=1}^{3}(1-ze^{i\psi_j})}.
\eeq
Similarly, $g'_n(z)$ converges uniformly on the compact subsets of the unit disk to
$$
g'(z)=ze^{i\theta}h'(z).
$$
If we take the logarithmic derivative of (\ref{BP1Ceq2}), we see that
$$
1+z\frac{h''(z)}{h'(z)}=\sum_{j=1}^3\frac{ze^{i\psi_j}}{1-ze^{i\psi_j}}+1, \quad z\in \ID.
$$
Since $w(z)=z/(1-z)$ maps $\ID$ onto the half plane  ${\rm Re\,}w>-1/2$, the last formula clearly
implies that
$${\rm Re}\left(1+z\frac{h''(z)}{h'(z)}\right)>-\frac{1}{2},  \quad z\in\ID.
$$
According to Theorem \ref{BP1Cth3}, $f$ is univalent close-to-convex in $\ID$.
This provides an alternate proof of Theorem 3.1 of Suffridge et. al. \cite{Suffridge-2001}.

\section{Applications of Lemmas \ref{BP1Ilem1} and \ref{BPlIlem2}}\label{BP1-sec3}
Consider the Gaussian hypergeometric function
\be\label{BP1eq7a}
 {}_2F_1(a,b;c;z):=F(a,b;c;z)=
\sum_{n=0}^{\infty}A_nz^n ,
\ee
where
$$ A_n=\frac{(a,n)(b,n)}{(c,n)(1,n)}.
$$
Here $a,b,c$ are complex numbers such that $c\neq -m$, $m=0,1,2,3, \ldots$, $(a,0)=1$ for $a\neq 0$ and,
for each positive integer $n$, $(a,n):=a(a+1)\cdots (a+n-1)$, see
for instance the recent book of Temme \cite{Temme:Special functions} and Anderson et. al. \cite{AVV}.
We see that $(a,n)=\Gamma(a+n)/\Gamma(a)$.
Often the Pochhammer notation $(a)_n$ is used instead of $(a,n)$.
In the exceptional case $c=-m$, $m=0,1,2,3, \ldots$, the function
$F(a,b;c;z)$ is clearly defined even if $a=-j$ or $b=-j$, where $j=0,1,2, \ldots$ and $j\leq m$.
The following well-known Gauss formula
\cite{Temme:Special functions} is crucial in the proof of our results of this section:
\begin{equation}\label{BP1eq8}
F(a,b;c;1) = \frac{\Gamma (c)\Gamma (c-a-b)}
{\Gamma (c-a)\Gamma (c-b)}< \infty \quad
\mbox{for ${\rm Re}\, c>{\rm Re}\,(a+b)$}.
\end{equation}

In order to generate nice examples of (fully) starlike and close-to-convex harmonic mappings, we consider
mappings whose co-analytic part involves Gaussian hypergeometric function.

\bthm\label{BP1th2}
Let either $a,b \in (-1, \infty)$ with $ab>0$, or $a,b\in\IC\setminus\{0\}$ with $b=\overline{a}$.
Assume that $c$ is a positive real number such that  $c>{\rm Re}\,(a+b)+1$, $\alpha\in\mathbb{C}$
and let
$$ f_k(z)=z+\overline {\alpha z^k F(a,b;c;z)} ~\mbox{ for $k=1,2$}.
$$
\begin{itemize}
\item[{\rm \textbf{(a)}}] If
\begin{equation}\label{BP1eq9}
\frac{\Gamma (c)\Gamma (c-a-b-1)}
{\Gamma (c-a)\Gamma (c-b)}[ab+2(c-a-b-1)]\leq \frac{1}{|\alpha|},
\end{equation}
where  $0<|\alpha|<1/2$, then $f_2\in \mathcal{S}_H^{*0}\cap {\mathcal C}_H^1$.

\item[{\rm \textbf{(b)}}] If
 \begin{equation}\label{BP1eq13}
\frac{\Gamma (c)\Gamma (c-a-b-1)}
{\Gamma (c-a)\Gamma (c-b)}[ab+c-a-b-1]\leq \frac{1+|\alpha|}{|\alpha|},
\end{equation}
where $2|\alpha|ab\leq c$, then
$f(z)=z+\overline {\alpha z (F(a,b;c;z)-1)} \in \mathcal{S}_H^{*0}\cap {\mathcal C}_H^1.$

\item[{\rm \textbf{(c)}}] If
$$\frac{\Gamma (c)\Gamma (c-a-b-1)}
{\Gamma (c-a)\Gamma (c-b)}[ab+c-a-b-1]\leq \frac{1}{|\alpha|},
$$
where $0<|\alpha|<1$, then $f_1\in {\mathcal C}_H^1$.
\end{itemize}
\ethm
\bpf
We present a proof of \textbf{(a)} and since the proofs of other two cases follow in a similar fashion,
we only include a key step for \textbf{(b)}.

\textbf{(a)} Set $h(z)=z$ and $\ds g(z)=\sum_{n=2}^{\infty}b_nz^n=\alpha z^2 F(a,b;c;z)$ so that
$$f_2(z)=z+\overline {\alpha z^2 F(a,b;c;z)}.
$$
By (\ref{BP1eq7a}) we have
\begin{equation}\label{BP1eq10}
b_n= \alpha A_{n-2}=\alpha\frac{(a,n-2)(b,n-2)}{(c,n-2)(1,n-2)} ~\mbox{ for }~ n\geq2.
\ee
By Lemma \ref{BPlIlem2} it suffices to show that
$K:=\ds \sum_{n=2}^{\infty}n|b_n|\leq1$. Using (\ref{BP1eq10}) it follows that
\begin{align*}
K &= |\alpha|\sum_{n=2}^{\infty}\frac{n(a,n-2)(b,n-2)}{(c,n-2)(1,n-2)}\\
&= |\alpha| \left( 2+ \sum_{n=1}^{\infty}\frac{(n+2)(a,n)(b,n)}{(c,n)(1,n)} \right)\\
&= |\alpha| \left(\frac{ab}{c} \sum_{n=1}^{\infty}\frac{(a+1,n-1)(b+1,n-1)}{(c+1,n-1)(1,n-1)}
  + 2\sum_{n=0}^{\infty}\frac{(a,n)(b,n)}{(c,n)(1,n)}  \right).
\end{align*}
By the hypothesis we have $c>a+b+1$ and both the series in the last expression converge and
so using the formula (\ref{BP1eq8}), we get
\begin{align*}
K&=|\alpha| \left(\frac{ab}{c} \left[\frac{\Gamma (c+1)\Gamma (c-a-b-1)}
{\Gamma (c-a)\Gamma (c-b)}\right]+ 2\frac{\Gamma(c) \Gamma(c-a-b)}{\Gamma(c-a) \Gamma(c-b)}
\right)\\
&=|\alpha| \left( \frac{\Gamma(c) \Gamma(c-a-b-1)}{\Gamma(c-a) \Gamma(c-b)}[ab+2(c-a-b-1)] \right).
\end{align*}
Clearly, (\ref{BP1eq9}) is equivalent to $K\leq1$. Thus, $f_2\in {\mathcal C}_H^1$ and is also
fully starlike in $\ID$. We complete the proof of \textbf{(a)}.

\textbf{(b)} For the proof of \textbf{(b)}, we consider $g$ defined by
$$g(z)=\alpha z (F(a,b;c;z)-1) = \alpha \sum_{n=2}^{\infty} A_{n-1}z^n
$$
and suffices to observe that
\begin{align*}
\alpha \sum_{n=2}^{\infty}n| A_{n-1}|
&=|\alpha| \left( \frac{ab}{c} \left[\frac{\Gamma (c+1)\Gamma
(c-a-b-1)} {\Gamma (c-a)\Gamma (c-b)}\right]+  \frac{\Gamma(c)
\Gamma(c-a-b)}{\Gamma(c-a) \Gamma(c-b)}-1 \right).
\end{align*}
\epf

The case $a=1$ of Theorem \ref{BP1th2}\textbf{(a)} and \textbf{(b)} gives

\bcor\label{BP1cor2}
Let $b$ and $c$ be positive real numbers and $\alpha$ be a complex number.
\begin{itemize}
\item[{\rm \textbf{(a)}}] If  $0<|\alpha|<1/2$ and
\begin{equation}\label{BP1eq11}
c\geq \beta^+ =\frac{3-6|\alpha|+(2-|\alpha|)b+\sqrt{|\alpha|^2(b^2+4b+4)+1+|\alpha|(4b^2-2b-4)}},
{2(1-2|\alpha|)}
\end{equation}
then
$$f_2(z)=z+\overline {\alpha z^2 F(1,b;c;z)} \in  \es_H^{*0}\cap {\mathcal C}_H^1.
$$
\item[{\rm \textbf{(b)}}] If $2|\alpha|b\leq c$ and
\begin{equation}\label{BP1eq14}
c\geq r^+=\frac{3+2b(1+|\alpha|)+\sqrt{b^2(4|\alpha|^2+4|\alpha|)+1}}{2},
\end{equation}
then $f(z)=z+\overline {\alpha z (F(1,b;c;z)-1)} \in \mathcal{S}_H^{*0}\cap {\mathcal C}_H^1$.
\end{itemize}
\ecor
\bpf
\textbf{(a)} Let $f_2(z)=z+\overline {\alpha z^2 F(1,b;c;z)}$.
It suffices to prove that if $c\geq \beta^+$
then the inequality (\ref{BP1eq9}) is satisfied with $a=1$.

It can be easily seen that $\beta^+>b+2$ and so the condition
$c\geq\beta^+$ implies that $c>b+2.$ Next, the condition (\ref{BP1eq9}) for $a=1$ reduces to
$$\frac{\Gamma (c)\Gamma (c-b-2)}
{\Gamma (c-1)\Gamma (c-b)}[b+2(c-b-2)]\leq \frac{1}{|\alpha|}
$$
which is equivalent to
$$(1-2|\alpha|)c^2+c[b(|\alpha|-2)-3+6\alpha]+b^2+(3-|\alpha|)b+2-4|\alpha|\geq0.
$$
Simplifying this inequality gives
$$
(1-2|\alpha|)(c-\beta^-)(c-\beta^+)\geq 0,
$$
where $\beta^+$ is given by (\ref{BP1eq11}) and
$$\beta^-=\frac{3-6|\alpha|+(2-|\alpha|)b-\sqrt{|\alpha|^2(b^2+4b+4)+1+|\alpha|(4b^2-2b-4)}}
{2(1-2|\alpha|)}.
$$
Since $\beta^+\geq\beta^-$ and by hypothesis $c\geq\beta^+,$ the inequality (\ref{BP1eq9}) holds.
It follows from Theorem \ref{BP1th2}\textbf{(a)} that $f_2\in{\mathcal C}_H^1$ and $f_2$
is fully starlike.

The proof for case \textbf{(b)} follows if one adopts a similar approach. In fact if we
set $a=1$ in Theorem \ref{BP1th2}\textbf{(b)}, then it is easy to see that the inequality
(\ref{BP1eq13}) is equivalent to
$$c^2+c(-2b(1+|\alpha|)-3)+(b^2+3b)(1+|\alpha|)+2 =(c-r^-)(c-r^+) \geq0
$$
where $r^+$ is given by (\ref{BP1eq14}) and
$$r^-=\frac{3+2b(1+|\alpha|)-\sqrt{b^2(4|\alpha|^2+4|\alpha|)+1}}{2}.
$$
Since $r^+\geq r^-$, the hypothesis that $c\geq r^+$ gives the desired conclusion.
\epf

As pointed out in Section \ref{sec-harmpol},  except the work of \cite{Suffridge-1998, Suffridge-2001},
there does not seem to exist a good technique to generate univalent harmonic polynomials. In view of
Theorem \ref{BP1th2}, we can obtain harmonic univalent polynomials that are close-to-convex and
fully starlike in $\ID$.

\bcor\label{BP1A2cor2}
Let $m$ be a positive integer, $c$ be a positive real number, $\alpha\in\mathbb{C}$, and let
$$ f_k(z)=z+\overline {\alpha z^k\sum_{n=0}^{m}\binom{m}{n}\frac{(m-n+1,n)}{(c,n)}z^n} ~\mbox{ for $k=1,2$.}
$$
\begin{itemize}
\item[{\rm \textbf{(a)}}]
If $0<|\alpha|<1/2$ and
$$
\frac{\Gamma (c)\Gamma (c+2m-1)}
{(\Gamma (c+m))^2}[m^2+2(c+2m-1)]\leq \frac{1}{|\alpha|},
$$
then $f_2\in \es_H^{*0}\cap {\mathcal C}_H^1$, and $f_2$ is indeed fully starlike in $\ID$.

\item[{\rm \textbf{(b)}}]
If $0<|\alpha|<1$ and
$$
\frac{\Gamma (c)\Gamma (c+2m-1)}
{(\Gamma (c+m))^2}[m^2+c+2m-1]\leq \frac{1}{|\alpha|},
$$
then $f_1\in {\mathcal C}_H^1$.
\end{itemize}
\ecor
\bpf
The results follow if we set $a=b=-m$ in Theorem \ref{BP1th2}\textbf{(a)},  and \textbf{(c)}, respectively.
\epf

On the other hand, Theorem \ref{BP1th2}\textbf{(b)} for $a=b=-m$ shows that
if $m$ is a positive integer, $c$ is a positive real number and  $\alpha\in\mathbb{C}$ is
such that $2|\alpha|m^2\leq c$ and
$$ 
\frac{\Gamma (c)\Gamma (c+2m-1)}
{(\Gamma (c+m))^2}[m^2+c+2m-1]\leq \frac{1+|\alpha|}{|\alpha|},
$$
then
$$ f(z)=z+\overline {\alpha z\sum_{n=1}^{m}\binom{m}{n}\frac{(m-n+1,n)}{(c,n)}z^n}
$$
belongs to $\es_H^{*0}\cap {\mathcal C}_H^1$, and $f$ is indeed fully starlike in $\ID$.



\beg\label{BP1A2eg2b}
If we let $m=3$ in Corollary \ref{BP1A2cor2}\textbf{(a)}, then we have the following: if $c$ is a
positive real number such that
$|\alpha|g(c)\leq 1,
$
where $\alpha$ is a complex number with $0<|\alpha|<1/2$ and
$$\ds g(c)=2+\frac{27}{c}+\frac{72}{c(c+1)}+\frac{30}{c(c+1)(c+2)},
$$
then the harmonic function
$$\ds f(z)=z+\overline {\alpha \left(z^2+\frac{9}{c}z^3+\frac{18}{c(c+1)}z^4
+\frac{6}{c(c+1)(c+2)}z^5\right)}
$$
is fully starlike in $\ID$.

Similarly, we see that if
$$
c\geq\frac{14|\alpha|-1+\sqrt{36|\alpha|^2+52|\alpha|+1}}{2(1-2|\alpha|)},
$$
where $\alpha$ is a complex number with $0<|\alpha|<1/2$, then the harmonic function
$$\ds f(z)=z+\overline {\alpha \left(z^2+\frac{4}{c}z^3+\frac{2}{c(c+1)}z^4\right)}
$$
belongs to $\es_H^{*0}\cap {\mathcal C}_H^1$ and is indeed fully starlike in $\ID$.
\eeg

\beg\label{BP1A2eg2d}
The choice $m=2$ in Corollary \ref{BP1A2cor2}\textbf{(b)} easily gives the following:
if $c$ is a positive real number such that
$$
c\geq\frac{9|\alpha|-1+\sqrt{25|\alpha|^2+38|\alpha|+1}}{2(1-|\alpha|)},
$$
where $\alpha$ is a complex number with $0<|\alpha|<1$, then
$$\ds f(z)=z+\overline {\alpha \left(z+\frac{4}{c}z^2+\frac{2}{c(c+1)}z^3\right)}
\in \mathcal{C}_H^1.
$$
\eeg

\bthm\label{BP1th5}
Let $a,b \in (-1, \infty) $. Assume that $c$ is a positive real number, $\alpha\in\mathbb{C}$ and,
for $k=0,1$, define
$$ f_k(z)=z+\overline{\alpha z^k\int\limits_{0}^{z}F(a,b;c;t)\,dt} .
$$
\begin{itemize}
\item[{\rm \textbf{(a)}}]
Let $ab>0$ or $a,b\in\IC\setminus\{0\}$ with $b=\overline{a}$, where $c>{\rm Re}\,(a+b)$ and $0<|\alpha|<1$
such that
$$
\frac{\Gamma (c)\Gamma (c-a-b)}
{\Gamma (c-a)\Gamma (c-b)}\leq \frac{1}{|\alpha|}.
$$
Then $ f_0\in\mathcal{C}_H^1.$
\item[{\rm \textbf{(b)}}]
Let $(a-1)(b-1)>0$,   or
$a,b\in\IC\setminus\{0,1\}$ with $b=\overline{a}$, where $c > \max\{1, {\rm Re}\,(a+b)\}$ and
$0<|\alpha|<1/2$ such that,
$$
  |\alpha|\frac{\Gamma(c)\Gamma(c-a-b)}{\Gamma(c-a)\Gamma(c-b)}
\left(1+\frac{c-a-b}{(a-1)(b-1)}\right)
    \leq1+|\alpha|\frac{(c-1)}{(a-1)(b-1)}.
$$
Then $ f_1 \in {\mathcal C}_H^1\cap \es_H^{*0}. $
Moreover, $f_1$ is fully starlike in $\ID$.

\end{itemize}
\ethm
\bpf
We give the proof of \textbf{(a)} and since the proof of \textbf{(b)} follows in a similar fashion,
we omit the details.

\textbf{(a)} Set $f_0(z):=z+\overline{g(z)}$, where
$$g(z)=\alpha\int\limits_{0}^{z}F(a,b;c;t)\, dt=\sum_{n=1}^{\infty}b_nz^n, \quad
b_n=\alpha\frac{(a,n-1)(b,n-1)}{(c,n-1)(1,n-1)n} ~\mbox{ for }~ n\geq 1.
$$
Therefore
\begin{align*}
\sum_{n=1}^{\infty}n|b_n|&=|\alpha|\sum_{n=1}^{\infty}\frac{(a,n-1)(b,n-1)}{(c,n-1)(1,n-1)}
=|\alpha| \frac{\Gamma(c)\Gamma(c-a-b)}{\Gamma(c-a) \Gamma(c-b)},
\end{align*}
as $c>{\rm Re}\,(a+b) $. The conclusion now follows from Lemma \ref{BPlIlem2}.
\epf

For instance, the case $a=1$ in Theorem \ref{BP1th5}{\rm \textbf{(a)}} shows that if
 $b$ and $c$ are positive real numbers such that
$$c\geq\frac{b+1-|\alpha|}{1-|\alpha|},
$$
where $\alpha$ is a complex number satisfying $0<|\alpha|<1$, then
$$ f(z)=z+\overline{\alpha\int\limits_{0}^{z}F(1,b;c;t)\,dt}\in\mathcal{C}_H^1.
$$

\bcor\label{BP1cor6}
Assume that $c$ is a positive real number and $\alpha$ is a complex number.
\begin{itemize}
\item[{\rm \textbf{(a)}}] Suppose that  either $a,b \in (-1, \infty)$ with $ab>0$, or
$a,b\in\IC\setminus\{0\}$ with $b=\overline{a}$. If
$c>{\rm Re}\,(a+b)+1$ and  $0<|\alpha|<1$ such that
\begin{equation}\label{BP1eq18a}
\frac{\Gamma (c+1)\Gamma (c-a-b-1)}
{\Gamma (c-a)\Gamma (c-b)}\leq \frac{1}{|\alpha|},
\end{equation}
then
$$ f(z)=z+\overline{(\alpha c/ (ab))[F(a,b;c;z)-1]}\in\mathcal{C}_H^1.
$$
\item[{\rm \textbf{(b)}}]  Suppose that  either $a,b\in(-1, \infty)$ with $ab>0$, or
 $a,b\in\IC\setminus\{0,1\}$ with $b=\overline{a}$. If $c > \max\{1, {\rm Re}\,(a+b)+1\}$ and
 $0<|\alpha|<1/2$ such that
$$
 |\alpha|\frac{\Gamma(c+1)\Gamma(c-a-b-1)}{\Gamma(c-a)\Gamma(c-b)}
\left(1+\frac{c-a-b-1}{ab}\right) \leq1+|\alpha|\frac{c}{ab},
$$
then
$$f(z)=z+\overline{[\alpha c/(ab)]z(F(a,b;c;z)-1)} \in {\mathcal C}_H^1\cap \es_H^{*0}.
$$
Moreover, $f$ is fully starlike in $\ID$.
\end{itemize}
\ecor
\bpf \textbf{(a)} The proof follows as a consequence of the following simple
identity for the first derivative of the hypergeometric function
$$
 abF(a+1,b+1;c+1;z)=cF'(a,b;c;z).
$$
Since
$$\int\limits_{0}^{z}F(a+1,b+1;c+1;t)\, dt=\frac{c}{ab}(F(a,b;c;z)-1),
$$
the conclusion follows if we apply Theorem \ref{BP1th5}{\rm \textbf{(a)}} and replace $a,b,c$ by
$a+1,b+1,c+1$, respectively.

The proof of case \textbf{(b)} follows if we apply Theorem \ref{BP1th5}{\rm \textbf{(b)}} with $a+1$, $b+1$, and $c+1$ instead
of $a$, $b$, and $c$, respectively.
\epf

\bcor\label{BP1cor7}
Let $\alpha$ be a complex number such that $0<|\alpha|<1$,
$b$ and $c$ be positive real numbers such that
\begin{equation}\label{BP1eq22}
c\geq r_1=\frac{3+2b-|\alpha|+\sqrt{|\alpha|^2(4b+1)+|\alpha|(8b+2)+1}}{2(1-|\alpha|)}.
\end{equation}
Then  $f(z)=z+\overline{(\alpha c/b)[F(1,b;c;z)-1]}\in\mathcal{C}_H^1$.
\ecor
\bpf
Let $f(z)=z+\overline{(\alpha c/b)[F(1,b;c;z)-1]}$. It suffices to prove that if $c\geq r_1$ then
inequality (\ref{BP1eq18a}) is satisfied with $a=1$. It is easily seen that $r_1>b+2$ and
hence, $c>b+2$ holds. Now the inequality (\ref{BP1eq18a}), with $a=1$ and
a simplification, is equivalent to
\begin{equation}\label{BP1eq23}
(1-|\alpha|)(c-r_1)(c-r_2)\geq 0,
\end{equation}
where $r_1$ is given by (\ref{BP1eq22}) and
$$r_2=\frac{3+2b-|\alpha|-\sqrt{|\alpha|^2(4b+1)+|\alpha|(8b+2)+1}}{2(1-|\alpha|)}.
$$
Since $r_1\geq r_2$ and by hypothesis $c\geq r_1$, the inequality (\ref{BP1eq23}) holds and thus,
by Corollary \ref{BP1cor6}\textbf{(a)}, $f$ belongs to $\mathcal{C}_H^1$,
and hence $f$ is close-to-convex in $\mathbb{D}.$
\epf


\bcor\label{BP1A2cor7}
Let $m$ be a positive integer, $c$ be a positive real number and $\alpha\in\mathbb{C}$. For $k\in \{0,1\}$, let
$$ f_k(z)=z+\overline {\alpha z^{k}\sum_{n=0}^{m}\binom{m}{n}\frac{(m-n+1,n)}{(c,n)}\frac{z^{n+1}}{n+1}}
$$
\begin{itemize}
\item[{\rm \textbf{(a)}}]
If $0<|\alpha|<1$ and
$$
\frac{\Gamma (c)\Gamma (c+2m)}
{(\Gamma (c+m))^2}\leq \frac{1}{|\alpha|},
$$
then $ f_0\in \mathcal{C}_H^1.$
\item[{\rm \textbf{(b)}}]If $0<|\alpha|<1/2$ and
$$
|\alpha|\frac{\Gamma(c)\Gamma(c+2m)}{(\Gamma(c+m))^2}
\left(1+\frac{c+2m}{(m+1)^2}\right)
    \leq1+|\alpha|\frac{(c-1)}{(m+1)^2},
$$
then $ f_1 \in {\mathcal C}_H^1\cap \es_H^{*0}.$
\end{itemize}
\ecor
\bpf Set $a=b=-m$ in Theorem \ref{BP1th5} {\rm \textbf{(a)}} and {\rm \textbf{(b)}}, respectively.
\epf

\beg\label{BP1A2eg2f}
Corollary \ref{BP1A2cor7}{\rm \textbf{(b)}} for $m=2$ gives the following:
if $c$ is a positive real number such that
$$
c\geq\frac{24|\alpha|-3+\sqrt{-48|\alpha|^2+168|\alpha|+9}}{6(1-2|\alpha|)},
$$
where $\alpha$ is a complex number with $0<|\alpha|<1/2$, then
$$\ds f(z)=z+\overline {\alpha \left(z^2+\frac{2}{c}z^3+\frac{2}{c(c+1)}\frac{z^4}{3}\right)}
\in \mathcal{C}_H^1\cap \es_H^{*0}.
$$.
\eeg

\subsection*{Acknowledgments}
The research of S. V. Bharanedhar was supported by the fellowship of the
University Grants Commission, India. The authors thank the referee for valuable comments
which improve the presentation of the results.

\end{document}